\documentstyle{czjphys209}

\makeatletter

\let\=\m@th

\newfont{\bslx}{cmbxsl10}

\newfont{\msamx}{msam10}
\newfont{\msamix}{msam9}
\newfont{\msamviii}{msam8}
\newfont{\msamvii}{msam7}
\newfont{\msamvi}{msam6}
\newfont{\msamv}{msam5}

\newfont{\msbmx}{msbm10}
\newfont{\msbmix}{msbm9}
\newfont{\msbmviii}{msbm8}
\newfont{\msbmvii}{msbm7}
\newfont{\msbmvi}{msbm6}
\newfont{\msbmv}{msbm5}

\newfont{\eufmx}{eufm10}
\newfont{\eufmix}{eufm9}
\newfont{\eufmviii}{eufm8}
\newfont{\eufmvii}{eufm7}
\newfont{\eufmvi}{eufm6}
\newfont{\eufmv}{eufm5}

\def\newsymbol#1#2#3#4#5{\let\next@\relax
 \ifnum#2=\@ne\let\next@\msafam@\else
 \ifnum#2=\tw@\let\next@\msbfam@\fi\fi
 \mathchardef#1="#3\next@#4#5}
\def\mathhexbox@#1#2#3{\relax
 \ifmmode\mathpalette{}{\m@th\mathchar"#1#2#3}%
 \else\leavevmode\hbox{$\m@th\mathchar"#1#2#3$}\fi}
\def\hexnumber@#1{\ifcase#1 0\or 1\or 2\or 3\or 4\or 5\or 6\or 7\or 8\or
 9\or A\or B\or C\or D\or E\or F\fi}

\newfam\frakfam \def\frak#1{{\fam\frakfam\relax#1}} \textfont\frakfam\eufmx
\scriptfont\frakfam\eufmvii \scriptscriptfont\frakfam\eufmv

\newfam\msafam \textfont\msafam\msamx
\scriptfont\msafam\msamvii \scriptscriptfont\msafam\msamv
\edef\msafam@{\hexnumber@\msafam}

\newfam\msbfam \textfont\msbfam\msbmx
\scriptfont\msbfam\msbmvii \scriptscriptfont\msbfam\msbmv
\edef\msbfam@{\hexnumber@\msbfam}

\def\Bbb#1{{\fam\msbfam\relax#1}}
\def\widehat#1{\setbox\z@\hbox{$\m@th#1$}%
 \ifdim\wd\z@>\tw@ em\mathaccent"0\msbfam@5B{#1}%
 \else\mathaccent"0362{#1}\fi}
\def\widetilde#1{\setbox\z@\hbox{$\m@th#1$}%
 \ifdim\wd\z@>\tw@ em\mathaccent"0\msbfam@5D{#1}%
 \else\mathaccent"0365{#1}\fi}

\newsymbol\leqslant 1336
\newsymbol\geqslant 133E

\def\roman#1{{\rm #1}}

\let\bls\baselineskip \let\nt\noindent \let\ignore\ignorespaces
\let\Hat\widehat \let\Tilde\widetilde \let\dollar\$ \let\ampersand\&
\let\tsize\textstyle \let\dsize\displaystyle \let\ssize\scriptstyle
\let\sss\scriptscriptstyle
  
\let\cline\centerline \let\lline\leftline \let\rline\rightline
\let\vp\vphantom \let\hp\hphantom \let\nt\noindent

\def\hline#1{\hbox to \hsize{#1}} \def\qqq{\qquad\quad}

\def\;{\relax\ifmmode\mskip\thickmuskip\else\kern.27777em\fi}
\def\!{\relax\ifmmode\mskip-\thinmuskip\else\negthinspace\fi}
\def\>{{\;\!}} \def\){{\;\;\!\!\!}} \def\]{{\;\!\!}} \def\}{{\]\]}}
\def\~{\raise.15ex\mhbox{-}}

\let\Sum\sum \def\sum{\Sum\limits}
\let\Prod\prod \def\prod{\Prod\limits}
\let\Int\int \def\int{\Int\limits}

\def\tsum{\mathop{\tsize\Sum}\limits}
\def\tprod{\mathop{\tsize\Prod}\limits}

\def\vsk#1>{\vskip#1\bls} 
\def\vv#1>{\vadjust{\vsk#1>}\ignore}
\def\vvn#1>{\vadjust{\nobreak\vsk#1>\nobreak}\ignore}
\def\vvgood{\vadjust{\penalty-500}}

\def\^#1{\hbox{\m@th#1}}
\def\h@ph{\discretionary{}{}{-}}
\def\$#1$-{\,\^{$#1$}\h@ph}

\def\mhbox#1{\hbox{\m@th$#1$}}

\def\gad#1{\global\advance#1\@ne}

\def\textindent#1{\indent\llap{#1\enspace}\ignorespaces}

\newcount\itemlet
\def\newbi{\itemlet 96} \newbi
\def\bitem{\gad\itemlet \par\hangindent1.5\parindent
 \hglue-.5\parindent\textindent{\rm \rlap{\char\the\itemlet}\hp{b})}}
\def\atem{\newbi\bitem}

\newcount\itemrm

\def\iitem{\gad\itemrm \par\hangindent1.5\parindent\hglue-.5\parindent
 \textindent{\upshape\hp{v}\llap{\romannumeral\the\itemrm})}}

\def\em#1{{\it #1\/}} \def\emph#1{{\sl #1\/}}

\def\Line#1{\kern-.5\hsize\line{\m@th$\dsize#1$}\kern-.5\hsize}
\def\Lline#1{\kern-.5\hsize\lline{\m@th$\dsize#1$}\kern-.5\hsize}
\def\Cline#1{\kern-.5\hsize\cline{\m@th$\dsize#1$}\kern-.5\hsize}
\def\Rline#1{\kern-.5\hsize\rline{\m@th$\dsize#1$}\kern-.5\hsize}

\def\Ll@p#1{\llap{\m@th$#1$}} \def\Rl@p#1{\rlap{\m@th$#1$}}
 \def\Cl@p#1{\llap{\m@th$#1$\hss}}
\def\Llap#1{\mathchoice{\Ll@p{\dsize#1}}{\Ll@p{\tsize#1}}{\Ll@p{\ssize#1}}%
 {\Ll@p{\sss#1}}}
\def\Clap#1{\mathchoice{\Cl@p{\dsize#1}}{\Cl@p{\tsize#1}}{\Cl@p{\ssize#1}}%
 {\Cl@p{\sss#1}}}
\def\Rlap#1{\mathchoice{\Rl@p{\dsize#1}}{\Rl@p{\tsize#1}}{\Rl@p{\ssize#1}}%
 {\Rl@p{\sss#1}}}
\def\LRtph#1#2{\setbox\z@\hbox{#1}\dimen\z@\wd\z@\hbox{\hbox to\dimen\z@{#2}}}
\def\LRph#1#2{\LRtph{\m@th$#1$}{\m@th$#2$}}

\def\Lto#1{\setbox\z@\mhbox{\tsize{#1}}%
 \mathrel{\mathop{\hbox to\wd\z@{\rightarrowfill}}\limits#1}}
\def\Lgets#1{\setbox\z@\mhbox{\tsize{#1}}%
 \mathrel{\mathop{\hbox to\wd\z@{\leftarrowfill}}\limits#1}}

\def\vpb#1{{\vp{\big(}}^{\]#1}} 

\let\alb\allowbreak 

 \let\x\times \let\ox\otimes 
\let\sub\subset  \let\tabs\+
\let\le\leqslant \let\ge\geqslant
 \let\8\infty \let\*\star
\let\bra\langle \let\ket\rangle
 
\let\map\mapsto  
 
 \def\vert{\ |\ } \def\nin{\not\in}

\let\lb\lbrace \let\rb\rbrace

\def\lsym#1{#1\alb\ldots\relax#1\alb}
\def\lc{\lsym,}   \def\lox{\lsym\ox}

\def\End{\mathop{\roman{End}\>}\nolimits}

\def\Res{\mathop{\roman{Res}\>}\limits}
\def\res{\mathop{\roman{res}\>}\limits}

\def\sgn{\mathop{\roman{sgn}\)}\limits}

\def\1{^{-1}} \let\underscore\_ \def\_#1{_{\Rlap{#1}}}
\def\vst#1{{\lower1.9\p@@\mhbox{\bigr|_{\raise.5\p@@\mhbox{\ssize#1}}}}}
\def\vrp#1:#2>{{\vrule height#1 depth#2 width\z@}}
\def\vru#1>{\vrp#1:\z@>} \def\vrd#1>{\vrp\z@:#1>}
\def\qqq{\qquad\quad} 
\def\sscr#1{\raise.3ex\mhbox{\sss#1}} \def\@@PS{\bold{OOPS!!!}}
 
\def\q{q\1} \def\2{^{-2}} \def\4{^{-4}}

\let\@ml@t\" \def\"#1{\ifmmode ^{(#1)}\else\@ml@t#1\fi}
\let\@c@t@\' \def\'#1{\ifmmode _{(#1)}\else\@c@t@#1\fi}
\let\colon\: \def\:{^{\vp{q}}} 

\let\al\alpha

 \let\Dl\Delta 
 \let\eps\varepsilon \let\epsilon\eps

\let\tht\theta \let\Tht\Theta

\let\ka\kappa
\let\la\lambda 

\let\si\sigma 
 
\let\pho\phi \let\phi\varphi

\def\C{\Bbb C}

\def\Z{\Bbb Z}

\def\Fc{{\cal F}}
\def\Fq{\Fc_{\!\!\]\sss e\]l\}l}}

\def\zh{\hat z}
\def\Fqh{\Rlap{\,\)\smash{\Hat{\!\]\phantom\Fc}}}\Fq}
\def\Fqt{\Rlap{\;\smash{\Tilde{\!\}\phantom\Fc}}}\Fq}
\def\FFhl{\Rlap\Fc\Rlap{\>\Hat{\phantom{\Fc}\;}}{\ssize\;}\Fq^{\)\ell}}
\def\Fqhl{\Fqh^{\ox\ell}}
\def\Mb{{\;\smash{\overline{\!\]\]M\!}\,}\vp M}}

\def\vb{\bar v}
\def\vti{\tilde v}
\def\wti{\tilde w}

\def\Ucup{{\tsize\bigcup}}

\def\Mx{M_{\rm ext}}

\def\slt{\frak{sl}_2} \def\Uq{U_q(\slt)}
\def\Cl{\C^{\)\ell}} \def\nums{\lb\)1\lc n\)\rb}

\def\zn{z_1\lc z_n} \def\zjn{z_1\lc p\)z_j\lc z_n}
\def\tell{t_1\lc t_\ell}
\def\Vox{V^{\ox n}} \def\sing{^{\sscr{\it sing}}}
\def\Vl{(\Vox)_\ell\:} \def\Vls{(\Vox)_\ell\sing}
\def\l@inf{\lower.21ex\mhbox{\ssize\8}} \def\9{_{\kern-.02em\l@inf}}
\def\Iell{I^{\ox\ell}} \def\Tq{T_{\]q^2}}
\def\0#1{\bra #1\>\ket}

\def\(#1){{\rm(\ref{#1})}}

\def\difl/{differential} \def\dif/{difference}
\def\cf.{cf.\ \ignore} \def\Cf.{Cf.\ \ignore}
\def\egv/{eigenvector} \def\eva/{eigenvalue} \def\eq/{equation}
\def\lhs/{the left hand side} \def\rhs/{the right hand side}
\def\Lhs/{The left hand side} \def\Rhs/{The right hand side}
\def\gby/{generated by} \def\wrt/{with respect to} \def\st/{such that}
\def\resp/{respectively} \def\off/{offdiagonal} \def\wt/{weight}
\def\pol/{polynomial} \def\rat/{rational} \def\tri/{trigonometric}
\def\fn/{function} \def\var/{variable} \def\raf/{\rat/ \fn/}
\def\inv/{invariant} \def\hol/{holomorphic} \def\hof/{\hol/ \fn/}
\def\mer/{meromorphic} \def\mef/{\mer/ \fn/} \def\mult/{multiplicity}
\def\sym/{symmetric} \def\perm/{permutation} \def\fd/{finite-dimensional}
\def\rep/{representation} \def\irr/{irreducible} \def\irrep/{\irr/ \rep/}
\def\hom/{homomorphism} \def\aut/{automorphism} \def\iso/{isomorphism}
\def\lex/{lexicographical} \def\as/{asymptotic} \def\asex/{\as/ expansion}
\def\ndeg/{nondegenerate} \def\neib/{neighbourhood} \def\deq/{\dif/ \eq/}
\def\hw/{highest \wt/} \def\gv/{generating vector} \def\eqv/{equivalent}
\def\msd/{method of steepest descend} \def\pd/{pairwise distinct}
\def\wlg/{without loss of generality} \def\Wlg/{Without loss of generality}
\def\onedim/{one-dimensional} \def\qcl/{quasiclassical} \def\hwv/{\hw/ vector}
\def\hgeom/{hyper\-geometric} \def\hint/{\hgeom/ integral}
\def\hwm/{\hw/ module} \def\emod/{evaluation module} \def\Vmod/{Verma module}
\def\symg/{\sym/ group} \def\sol/{solution} \def\eval/{evaluation}
\def\anf/{analytic \fn/} \def\anco/{analytic continuation}
\def\qg/{quantum group} \def\qaff/{quantum affine algebra}

\def\Rm/{\^{$R$-}matrix} \def\Rms/{\^{$R$-}matrices} \def\YB/{Yang-Baxter \eq/}
\def\Ba/{Bethe ansatz} \def\Bv/{Bethe vector} \def\Bae/{\Ba/ \eq/}
\def\KZv/{Knizh\-nik-Zamo\-lod\-chi\-kov} \def\KZvB/{\KZv/-Bernard}
\def\KZ/{{\sl KZ\/}} \def\qKZ/{{\sl qKZ\/}}
\def\KZB/{{\sl KZB\/}} \def\qKZB/{{\sl qKZB\/}}
\def\qKZo/{\qKZ/ operator} \def\qKZc/{\qKZ/ connection}
\def\KZe/{\KZ/ \eq/} \def\qKZe/{\qKZ/ \eq/} \def\qKZBe/{\qKZB/ \eq/}
\def\qKZb/{{\bslx q\)KZ\/}}

\def\h@ph{\discretionary{}{}{-}} \def\$#1$-{\,\^{$#1$}\h@ph}

\def\Asym{\mathop{\roman{Asym}\>}\nolimits}
\let\beq\be \let\eeq\ee \let\nn\nonumber
\def\kme{\kern-20pt}

\newtheorem{lm}{Lemma}[section]
\newtheorem{prop}[lm]{Proposition}
\newtheorem{th}[lm]{Theorem}

\let\secti@n\section \def\section{\setcounter{equation}0\secti@n}

\makeatother

\begin{document}

\title{HYPERGEOMETRIC SOLUTIONS OF THE \qKZb/ EQUATION\\[4pt]
AT LEVEL ZERO\footnote{Presented at the 8th Coloquium ``Quantum groups and
integrable systems'', Prague 17\)--\)19 June \)1999}}
\authori{V.\;Tarasov}
\addressi{St.\,Petersburg Branch of Steklov Mathematical Institute\\[2pt]
St.\,Petersburg \,191011, Russia}

\evensidemargin\oddsidemargin

\authorii{}
\addressii{}
\authoriii{}
\addressiii{}

\headtitle{Hypergeometric solutions of the \qKZ/ equation at level zero}
\headauthor{V.\;Tarasov}
\specialhead{V.\;Tarasov:
Hypergeometric solutions of the \qKZ/ equation at level zero}

\evidence{}
\daterec{\;18 August \)1999}	
\cislo{1}  \year{2000}
\setcounter{page}{193}
\pagesfromto{193--200}

\maketitle

\begin{abstract}
We discuss relations between different integral formulae for solutions
of the quantized Knizhnik-Zamolodchikov (\qKZ/) equation at level zero in
the $\Uq$ case for ${|\)q\)|<1}$. Smirnov type formulae of M.\,Jimbo et al.\
are derived from the general approach of A.\,Varchenko and the author.
The consideration is parallel to the \qKZ/ equation in the rational
$\slt$ case done by A.\,Nakayashiki, S.\,Pakuliak and the author.
\end{abstract}

\section{Introduction}
The quantized \KZv/ (\qKZ/) \eq/ is a holonomic system of \deq/s for
a \fn/ taking values in a tensor product of \rep/s of a Lie algebra or the
corresponding \qg/. In this note we consider the \qKZe/ associated with the
\qg/ $\Uq$, see \(qKZ). We concentrate on the case of the \qKZe/ at level zero,
which means a particular relation between $q$ and the step $p$ of the \eq/:
$p=q^4\!$. This type of the \qKZe/ was inroduced by F.\,Smirnov \cite{S} as
\eq/s for form factors in integrable models of quantum field theory.

We are interested in the \hgeom/ \sol/s of the \qKZe/. In the case in question
such \sol/s were given by Smirnov and Jimbo et al., see \cite{JKMO}.
The general construction of the \hgeom/ \sol/s has been developed by Varchenko
and the author \cite{TV}. It turns out that the general approach applied to
the level zero case gives formulae which are similar to Smirnov type formulae,
but do not coincide with them. On the other hand, all attempts to extend
Smirnov type formulae to the case of arbitrary level failed, which indicates an
intimate relation of these formulae to the special features of the level zero
case.

The aim of the note is to derive Smirnov type formulae for the \hgeom/ \sol/s,
\cf. \(detQ) and \(det0), from the general type formula \(PsiW). In the \rat/
$\slt$ case this has been done in \cite{NPT}. Here, we consider the case
${|\)p\)|<1}$. The case $|\)p\)|=1$ has been recently studied by Y.\,Takeyama.

Almost all proofs in the case in question are nearly the same as in the \rat/
case, see \cite{NPT}. We also refer a reader to this paper for more detailed
comments, motivations and explanations of the constructions. Completeness of
the \hgeom/ \sol/s in the \rat/ case at level zero has been proved in
\cite{T}\); in the case to be discussed it can be shown more or less similarly.

Formula \(PsiW) for the \hgeom/ solutions contains three main ingredients:
\raf/s $w_M$, an elliptic \fn/ $W$ and the \hint/ $\Iell$, which pairs
the \fn/s. We describe the constituents subsequently in Sections 3 and 4.
In Section 5 we describe the \hgeom/ \sol/s, \cf. Theorem \(hsol),
and derive Smirnov type formulae for \sol/s.

\section{The \qKZb/ equation at level zero}
Let $V=\C\>v_+\oplus\C\>v_-$. Introduce matrices
$$
\si^+=\}\left({0\ \;1\atop 0\ \;0}\right),\quad
\si^-=\}\left({0\ \;0\atop 1\ \;0}\right),\quad
\tau^+=\}\left({1\ \;0\atop 0\ \;0}\right),\quad
\tau^-=\}\left({0\ \;0\atop 0\ \;1}\right).
\vv.2>
$$
Fix a complex number $q\ne 0,1$. Let $R(z)\in\End(V^{\ox2})$ be the following
\Rm/:
\bea
R(z)\,=\,\tau^+\!\ox\tau^++\)\tau^-\!\ox\tau^-&\kme& {}+\,
{z-1\over qz-\q\!}\,(\tau^+\!\ox\tau^-+\tau^-\!\ox\tau^+)\>+{}\nn
\\[4pt]
&\kme& {}+\,{q-\q\over qz-\q\!}\,(z\>\si^+\!\ox\si^-+\si^-\!\ox\si^+)\,.\nn
\vv.1>
\eea
Fix complex numbers ${p\ne 0,1}$ and ${\ka\ne 0}$. We consider the \qKZe/ for
a \$\Vox$-valued \fn/ $\Psi(\zn)$:
\beq
\Psi(\zjn)\,=\,K_j(\zn)\>\Psi(\zn)\,,
\label{qKZ}
\vv-.6>
\eeq
\bea
K_j(\zn)\,&\kme&{}=\,R_{j,j-1}(p\)z_j/z_{j-1})\ldots R_{j,1}(p\)z_j/z_1)\,\x{}
\nn\\[4pt]
&\kme&{}\>\x\,
(\tau^+_j\}+\ka\)\tau^-_j)\>R_{j,n}(z_j/z_n)\ldots R_{j,j+1}(z_j/z_{j+1})\,,
\kern-2em\nn
\eea
$j=1\lc n$. The \em{level} of the \qKZe/ is determined by
$p=q^{\)2(level\)+\)2)}$.

\vsk.2>
For any $\ell$ \st/ $0\le\ell\le n$ denote by $\Vl$ the \wt/ subspace:
$$
\Vl\>=\>\bigl\lb\)v\in\Vox\vert\tsum_{j=1}^n\tau^-_jv=\ell\)v\)\bigr\rb\,.
\vv-.3>
$$
The \qKZe/ \(qKZ) respects the weight decomposition of $\Vox$.

\vsk.2>
Consider the \qg/ $\Uq$ with generators $e,f,k$ and the coproduct
$$
\Dl(e)\)=\)k\1\!\ox e+e\ox 1\,,\qquad
\Dl(f)\)=\)f\ox k+1\ox f\,,\qquad \Dl(k)\)=\)k\ox k\,.
$$
$\Uq$ acts in $V$ by the rule: $e\map\si^+\!$, \,$f\map\si^-\!$,
\vv.15>
\,$k\map q\)\tau^+\!+\q\tau^-\!$, \,and in $\Vox\!$ according to
the coproduct $\Dl$.
\vv.15>
If $\ka=q^{\)2\ell-2-n}\!$, then the \qKZe/ preserves the subspace
$\Vls\sub\Vl$ of \$\Uq$-singular vectors:
$$
\Vls=\>\bigl\lb\)v\in\Vl\vert e\)v=0\)\bigr\rb\,.
$$

In this note we consider the case of the \qKZe/ at level zero, i.e.\ all over
the paper we assume that $p=q^4$. Furthermore, fixing an integer $\ell$ \st/
$0\le\ell\le n$ we discuss \sol/s of the \qKZe/ \(qKZ) taking values in
\vv.1>
the \wt/ subspace $\Vl$. We will pay special attention to the case
\vvgood
$\ka=q^{\)2\ell-2-n}\!$. In this case we assume that $2\ell\le n$ and consider
\sol/s of the \qKZe/ taking values in the subspace $\Vls$ of singular vectors.

We assume that $|p|<1$. This is important for the analytic part of the story
in Sections 4 and 5.

\vsk->
\vsk0>
\section{Rational functions}
Let ${M=\lb\)m_1\lsym<m_\ell\)\rb}$ be a subset of $\nums$. The subset defines
a point $\zh_M=(z_{m_1}\lc z_{m_\ell})\in\Cl\!$. For subsets $M,N$
we say that $M\le N$ if $\#M=\#N$ and $m_a\le n_a$ for any $a=1\lc\#M$.

\vsk.2>
Given a \fn/ $f(\tell)$ we set
$$
\Asym f(\tell)=
\tsum_{\si\in\ibf S\_\ell}\sgn(\si)\)f(t_{\si_1}\lc t_{\si_\ell})\,,
$$
and for any point $u=(u_1\lc u_\ell)\in\Cl$ we define
$$
\Res f(u)\,=\,\res\bigl(\}\ldots\)\res\bigl(\)t_1\1\ldots t_\ell\1\]f(\tell)
\bigr)\]\big|_{t_\ell=u_\ell}\)\ldots\)\bigr)\]\big|_{t_1=u_1}\,.
$$
For any ${M\sub\nums}$ and ${m\in M}$ let $\mu\"m_M\!,\,g_M,\,w_M$ and $\wti_M$
be the following \fn/s:
$$
\mu\"m_M(t)\,=\,{t\over t-z_m}\>
\prod_{\tsize{j\in M\atop j\ne m}}\]{t-q^2z_j\over z_m-q^2z_j}\;,
\vv-.2>
$$
$$
\wti(\tell)\,=\Asym\Bigl(\,\prod_{a=1}^\ell\mu\"{m_a}(t_a)\]\Bigr)\,=\,
\det\}\bigl[\)\mu\"{m_a}_M(t_b)\)\bigr]_{a,b=1}^\ell\,,
$$
$$
g_M(\tell)\,=\,\prod_{a=1}^\ell\Bigl(\>{t_a\over t_a-z_{m_a}\!}\,\,
\prod_{1\le j<m\_a}{q\1t_a-qz_j\over t_a-z_j}\>\Bigr)
\prod_{1\le a<b\le\ell}(q\2 t_a-t_b)\,,
$$
\vvn.3>
$$
w_M\,=\,\Asym g_M\,.
$$
\begin{lm}
\label{Res}
Let ${M,N\sub\nums}$, ${\#M=\#N}$. Then
$$
\Res\wti_M(\zh_M)\>=\>1\,,\qqq \Res w_M(\zh_M)=\Res g_M(\zh_M)\,,
$$
${\Res\wti_M(\zh_N)=0}$ for ${M\ne N}$ \,and \ ${\Res w_M(\zh_N)=0}$ unless
$N\le M$. Moreover
$$
\tsize w_M\,=\sum_{N\le M}\}\wti_M\>\Res w_M(\zh_N)\,.
$$
\end{lm}
\begin{lm}
\label{D1}
For any $M\sub\nums$, $\#M=\ell-1$, the following relation holds\/{\rm:}
\bea
&\kme&
(q-\q)\sum_{k\nin M}^{}q^{\)2\la_M(k)-k}\>w_{M\)\cup\)\lb k\rb}(\tell)\,={}\nn
\\[3pt]
&\kme& {}=\,\Asym\biggl(\Bigl(\,\prod_{a=2}^\ell\)(q\2t_1-t_a)\,-
\prod_{a=2}^\ell\)(q^2t_1-t_a)\,\prod_{j=1}^n\>{q\2t_1-z_j\over t_1-z_j}
\>\Bigr)\>g(t_2\lc t_\ell)\biggr)\,.\nn
\eea
Here \,$\la_M(k)=\#\lb\)m\in M\vert m<k\)\rb$.
\end{lm}
Let $\Mx=\lb\)1\lc\ell\)\rb$. Say that $\Mx$ is the \em{extremal} subset.
Due to Lemma \ref{Res}
\vvn.2>
$$
w_{\Mx}(\tell)\,=\prod_{1\le a<b\le\ell}\!
{(q\2z_{m_a}\!-z_{m_b})\>(z_{m_a}\!-q\2z_{m_b})\over q\1(z_{m_a}\!-z_{m_b})}
\ \wti_{\Mx}(\tell)\,.
$$
\vsk-.4>
Let $\Fc$ be the space of \fn/s $f(t)$ \st/ the product
\vv-.15>
$f(t)\prod_{j=1}^n(t-z_j)$ is a Laurent \pol/ in $t$.
We denote by $\Fc^{\ox k}\!$ the space of \fn/s in $k$ \var/s \st/ for each
given \var/ these \fn/s considered as \fn/s of the distinguished \var/ belong
to $\Fc$. Other tensor products of spaces of \fn/s are to be understood
similarly. It is clear that for any subset $M$ of cardinality $\ell$
the \fn/s $g_M\,, w_M\,, \wti_M$ belong to $\Fc^{\ox\ell}\!$.

\vsk-.3>
\vsk0>
\section{Elliptic functions and the \hgeom/ integral}
\vsk-.3>

Let \,$(u)\9=\prod_{s=0}^\8(1-p^su)$ and let \,$\tht(u)=(u)\9\>(p/u)\9\>(p)\9$
\,be the Jacobi theta-\fn/. Denote by $\Fqt$ the space of \fn/s $F(t)$
\st/ the product $F(t)\prod_{j=1}^n\tht(t/z_j)$ is \hol/ for $t\ne 0$.

Let ${\al=\ka\>q^{n-2\ell+2}\!}$.
\vv.1>
Set ${\Fq=\lb\)F\in\Fqt\vert F(p\)t)=\al\)F(t)\)\rb}$. For even $n$ also set
${\Fqh=\lb\)F\in\Fqt\vert F(p\)t)-\al\)F(t)\in\C\>\Tht(t)\)\rb}$ where
\beq
\Tht(t)\,=\,t^{-n/2}\>\prod_{j=1}^n\,{\tht(q\2t/z_j)\over\tht(t/z_j)}\;.
\label{Tht}
\vv.1>
\eeq
It is easy to see that $\dim\)\Fq=n$ and $\dim\)\Fqh=n+1$. Notice that
$\Tht\in\Fqh$.

\vsk.5>
\nt
{\sl Remark.}\enspace
In the general situation one needs only an analogue of the space $\Fq$ of
quasiperiodic \fn/s, see \cite{TV}. However, in the case in question it turns
out that sometime the space $\Fq$ is not enough to produce all, or even
nonzero, \hgeom/ \sol/s of the \qKZe/ \(qKZ). This actually happens if
${\al=q^{\)2k}}$ for an integer $k$. The proposed extension $\Fqh$ is adapted
to the case $\al=1$.

\vsk.3>
Let \ ${\dsize\pho(t)\>=\)\prod_{j=1}^n\,{(t/z_j)\9\over(q\2t/z_j)\9\!}}$ \ be
\vv-.2>
the \em{phase \fn/}. Introduce the sets:
$\Ucup^+_s=\lb\)p^{\)k}q^{\)2}z_j\vert j=1\lc n\,,\ {k\in\Z_{\le s}}\)\rb$
\vv.1>
\,and \;${\Ucup^-_s=\lb\)p^{-k}z_j\vert j=1\lc n\,,\ k\in\Z_{\le s}\)\rb}$.
The \hint/ $I(f,F)$ is given by the formula
\vvn.2>
\beq
I(f,F)\,=\,{1\over 2\pi i}\,\int_{\}C\>}\pho\,f\>F\>{dt\over t}
\label{hint}
\eeq
where $C$ is a simple closed curve oriented counterclockwise and separating
the sets $\Ucup^+_1$ and $\Ucup^-_1$. For \fn/s in $\ell$ \var/s we write
\vvn.1>
\beq
\Iell(w,W)\,=\,{1\over(2\pi i)^\ell}\,
\int_{C^\ell}w(\tell)\,W(\tell)\>\tprod_{a=1}^\ell\pho(t_a)\,{dt_a\over t_a}\;.
\eeq
\begin{prop}
For any $f\in\Fc$ and $F\in\Fqt$ the \hint/ $I(f,F)$ is well defined and
does not depend on a particular choice of the contour $C$.
\end{prop}
\vsk->
\vsk0>
\begin{lm}
\label{I0}
Let ${f\in\Fc}\!$. If $f(t)=O(t)$ as $t\to 0$, then $I(f,1)=0$.
For even $n$, if $f(t)=O(t^k)$ as $t\to\8$, ${2k\le n}$, and $\Tht$ is
given by formula \(Tht), then $I(f,\Tht)=0$.
\end{lm}
\vsk-.5>
Denote by $D$ the operator defined by
\vvn-.5>
$$
Df(t)\,=\,f(t)-\al\>f(p\)t)\prod_{j=1}^n {q\2t-z_j\over t-z_j}\;.
\vv-.1>
$$
The \fn/s of the form $Df$ are called the \em{total \dif/s}.
\begin{prop}
\label{ID}
Let ${f\in D\Fc}$ or ${(f-Df)\in\Fc}$.
\vsk.15>
\atem
For any ${F\in\Fq}$ we have that $I(Df,F)=0$.
\bitem
Let n be even, and $f(t)=O(t^k)$ as $t\to\8$, ${2k\le n}$.
Then for any ${F\in\Fqh}$ we have that $I(Df,F)=0$.
\end{prop}

\vsk->
\vsk0>
\section{Hypergeometric solutions of the \qKZb/ equation}
For any subset ${M\sub\nums}$ define a vector ${v_M\in(\Vox)_{\#M}}$ by
the rule: $v_M=v_{\eps_1}\lox v_{\eps_n}$ \>where $\eps_j=+$ for $j\nin M$ and
\vv.1>
$\eps_j=-$ for $j\in M$. For any $W\in\Fqhl$ we set
\vv-.2>
\beq
\Psi_W\,=\sum_{\#M=\ell}\Iell(w_M,W)\,v_M\,.
\label{PsiW}
\vv-.1>
\eeq
It is clear that $\Psi_{\Asym W}=\ell\)!\,\Psi_W$. By Lemma \ref{I0} $\Psi_W=0$
if $W\!\in\Fqh^{\ox(\ell-1)}\ox\C\>1$, or $W\!\in\Fqh^{\ox(\ell-1)}\ox\C\>\Tht$,
$n$ is even and $2\ell\le n$.
\begin{lm}
Let $\al=1$. Then $\Psi_W\in\Vls$ for any $W\!\in\Fqhl$.
\end{lm}
\vsk-.3>
The statement follows from Lemma \ref{D1} and Proposition \ref{ID}.

\vsk.3>
Let $\FFhl$ be the space of \mef/s $W(\tell;\zn)$ with the properties:
\vsk.2>
\atem
$W\!\in\Fqhl$ as a \fn/ of $\tell$ for any given $\zn$;
\vsk.2>
\bitem
$W(\tell;\zjn)\>=\>q^{-\ell}\>W(\tell;\zn)$ \;for any $j$.

\begin{th}
\label{hsol}
For any $W\!\in\FFhl$ the \fn/ $\Psi_W(\zn)$ is a \sol/ of
\vv.1>
the \qKZe/ \(qKZ) with $p=q^4$ and $\ka=\al\>q^{\)2\ell-2-n}\!$
taking values in $\Vl$.
\end{th}

\vsk-.2>
$\Psi_W(\zn)$ is called a \em{\hgeom/ \sol/} of the \qKZe/.

\vsk.5>
Theorem \ref{hsol} follows from the results on formal integral \rep/s for
\sol/s of the \qKZe/ \cite{V} and Proposition \ref{ID}.

\goodbreak
\vsk.3>
Take another basis in $\Vl$: $\vti_M=\!\sum_{N\ge M}\!v_N\>\Res w_N(\zh_M)$,
\vvn-.5>
$\#M=\ell$. Then we have
\vvn.2>
$$
\Psi_W\,=\sum_{\#M=\ell}\Iell(\wti_M,W)\;\vti_M\,.
$$
For $W(\tell)=W_1(t_1)\ldots W_\ell(t_\ell)$, $W_1\lc W_\ell\in\Fqh$, the last
formula can be written in the determinant form:
\beq
\Psi_W\,=\sum_{\#M=\ell}
\det\]\bigl[I(\mu\"{m_a}_M,W_b)\bigr]_{a,b=1}^\ell\,\vti_M\,.
\label{det}
\eeq
\vsk-.3>
{}From now on until the end of the section we assume that $\al=1$, so that
$\ka=q^{\)2\ell-2-n}\!$. In this special case $\Psi_W$ can be also written
via suitable \pol/s rather than \raf/s. For any $M\sub\nums$ set
$$
P_M^+(t)\,=\prod_{m\in M}\!(q\4t-z_m)\,,\qqq
P_M^-(t)\,=\prod_{k\nin M}\!(q\4t-z_k)\,.
\vv-.2>
$$
Denote by $\Tq$ the following operator: $\Tq f(t)=f(t)-f(q^2t)$. For any \raf/
$f(t)$ let $[f(t)]_+$ be its \pol/ part. Define \pol/s $Q_M\"1\lc Q_M\"\ell$
by the rule
\vv-.7>
\bea
Q_M\"a(t)\, &\kme&{} =\,q^{\)4a}P_M^-(t)\>
\Bigl[\)\Tq\Bigl(\){P_M^+(t)\over t^a}\)\Bigr)\Bigr]_+ +{}\nn
\\[6pt]
&\kme& {}\)+\,q^{\)2a}P_M^+(q^2t)\)
\Bigl[\)\Tq\Bigl(\){P_M^-(t)\over P_M^+(q^2t)}\,
\Bigl[\>{P_M^+(q^2t)\over t^a}\>\Bigr]_+\)\Bigr)\Bigr]_+.
\label{QM}
\eea
The \pol/ $Q_M\"a(t)$ essentially coincides with the \pol/
\vv.2>
$A_{\ell+1-a}\"{\ell,n-\ell}(t\)|\zh_M|\zh_{\Mb})$ from \cite{JKMO}.
Here $\Mb=\nums\setminus\}M$ and $\zh_\Mb$is defined similarly to $\zh_M$.
\begin{prop}
For any ${M\sub\nums}$, ${\#M=\ell}$, and any ${m\in M}$ the following identity
holds\/{\rm:}
\vv-.2>
$$
D\Bigl(\prod_{\tsize{k=1\atop k\ne m}}^n(q\4t-z_k)\]\Bigr)\,=\,
-\,z_m\1\>\prod_{k=1}^n\](q\2z_m-z_k)\,\mu_M\"m(t)\;+\]
\sum_{a=1}^\ell\>Q_M\"a(t)\,z_m^{a-1}\>.
$$
\end{prop}
\vsk-.1>
Set
\vvn-.5>
$$
\vb_M\,=\!\prod_{1\le a<b\le\ell}{z_{m_a}\!-z_{m_b}\over
(q\2z_{m_a}\!-z_{m_b})\>(z_{m_a}\!-q\2z_{m_b})}\ \vti_M\,,
$$
which provides that $\,{\vb_{\Mx}\}-\)q^{\)\ell(\ell-2)/2}\>v_{\Mx}}$
is a linear combination of vectors $v_M$ with $M\ne\Mx$, $\#M=\ell$. Then
\vvgood
by the last Proposition and Lemma \ref{ID} the formula \(det) transforms to
\vvn.1>
\beq
\kern-.5em \Psi_W\,=\,(q\2\!-1)\vpb{\]-\ell}\!\sum_{\#M=\ell}\]
\prod_{\tsize{k\nin M\atop m\in M}}\!{1\over q\2z_m-z_k\}}
\ \det\]\bigl[I(Q_M\"a,W_b)\bigr]_{a,b=1}^\ell\;\vb_M \kern-.25em
\label{detQ}
\vvn-.1>
\eeq
provided that $W_1\lc W_\ell\in\Fq$. This formula matches Theorem 4.5
in \cite{JKMO}.

\vsk.3>
Observe that if $Q$ is a \pol/ and ${W\in\Fq}$, then the integrand of
the integral $I(Q,W)$ has no poles at points $z_j\,,\,p\1z_j\,,\,p\>q^2z_j$
for any ${j=1\lc n}$, and therefore, we have that
$$
I(Q,W)\,=\,{1\over 2\pi i}\,\int_{\}C'\)}\pho\,Q\>W\>{dt\over t}
$$
for any contour $C'$ separating the sets $\Ucup^+_0$ and $\Ucup^-_{-1}$.
\vv.15>
For instance, if $|z_1|\lsym=|z_n|=1$, then one can take $C'$ to be a circle
$|t|=q^3$.

\vsk.3>
{}From now on let $2\ell=n$, and recall that $\al=1$. For any $F\in\Fqh$
define its \em{discrepancy} $\0F$ by the rule: $\0F\>\Tht(t)=F(t)-F(p\)t)$.
For any subset $M$ set $\xi_M=DP_M^-$. There are two important relations
involving $\xi_M$:
$$
\xi_M(t)\,=\sum_{m\in M}\)\mu_M\"m(t)\,z_m\1\]\res\xi_M(z_m)
\vv-.2>
$$
and $I(\xi_M,F)=-\,q^{-4\ell}\0F$ for any $F\in\Fqh$. Using the former equality
to replace $\mu_M\"{m_\ell}$ by $\xi_M$ in \(det) and integrating the terms
with $\xi_M$ according to the latter one, we obtain
\vv-.7>
$$
\Psi_W\,=\,\0{W_\ell}\}\sum_{\#M=\ell}\}z_{m_\ell}\,{\prod_{a=1}^{\ell-1}
(z_{m_\ell}\!-z_{m_a})\over\prod_{k=1}^n(z_{m_\ell}\!-q^2z_k)}
\ \det\]\bigl[I(\mu\"{m_a}_M,W_b)\bigr]_{a,b=1}^{\ell-1}\,\vti_M
\vv.1>
$$
provided that $W_1\lc W_{\ell-1}\in\Fq$.
In the case $2\ell=n$ formula \(QM) for the \pol/s $Q_M\"a$ simplifies:
$$
Q_M\"a(t)\,=\,q^{\)4a}P_M^-(t)\>\bigl[\)\Tq\bigl(t^{-a}P_M^+(t)\bigr)\bigr]_+
+ q^{\)2a}P_M^+(q^2t)\>\bigl[\)\Tq\bigl(t^{-a}P_M^-(t)\bigr)\bigr]_+\,.
$$
Moreover, the \pol/ $Q_M\"\ell$ vanishes identically. Finally, we get
an analogue of formula \(detQ)\):
\vvn-.3>
\beq
\label{det0}
\kern-.6em
\Psi_W\,=\,(q^{\)2}\!-q^{\)4})\vpb{-\ell}\,\0{W_\ell}\}\sum_{\#M=\ell}\]
\prod_{\tsize{k\nin M\atop m\in M}}\!{1\over q\2z_m-z_k\}}
\ \det\]\bigl[I(Q_M\"a,W_b)\bigr]_{a,b=1}^{\ell-1}\;\vb_M\,.
\vvn-.2>
\eeq
matching the modification (4.17) to Theorem 4.5 in \cite{JKMO}.

\end{document}